\documentclass[reqno,12pt]{amsart}

\usepackage{fullpage}
\usepackage[french]{babel}
\usepackage[T1]{fontenc}
\usepackage{graphicx}
\usepackage{amsmath}
\usepackage{amsfonts}
\usepackage{amssymb}
\usepackage{a4wide}

\newtheorem{theo}{Théorème}
\newtheorem*{Crit1}{Critère pour les applications miroir}
\newtheorem*{Crit2}{Critère pour les applications de type miroir}

\newtheorem*{theo2}{Théorème}

\newtheorem{cor}{Corollaire}

\newtheorem{lemme}{Lemme}

\theoremstyle{remark}
\newtheorem*{Remarque}{Remarque}

\numberwithin{equation}{section}

\author{E. DELAYGUE}
\title{Intégralité des coefficients de Taylor de racines d'applications miroir}
\date{}

\begin{document}
\maketitle

\begin{abstract}
Nous démontrons l'intégralité des coefficients de Taylor de racines de séries de la forme $q(z):=z\exp(G(z)/F(z))$, où $F(z)$ et $G(z)+\log(z)F(z)$ sont des solutions particulières de certaines équations différentielles hypergéométriques généralisées. Cela nous permet de démontrer une conjecture de Zhou énoncée dans \og Integrality properties of variations of Mahler measures \fg [arXiv:1006.2428v1 math.AG]. La preuve de ces résultats est une adaptation des techniques utilisées dans notre article \og Critère pour l'intégralité des coefficients de Taylor des applications miroir \fg\,[J. Reine Angew. Math. (à paraître)].
\end{abstract}

\section{Introduction}

\subsection{\'Enoncés des résultats principaux}\label{section mirror map}

Nous continuons l'étude des coefficients de Taylor d'applications miroir d'origine hypergéométrique commencée dans \cite{Delaygue} où nous avons établi des critères d'intégralité pour les coefficients de Taylor d'applications miroir et d'applications de \textit{type miroir}. Nous allons raffiner notre critère pour les applications de type miroir, ce qui nous permettra d'obtenir, dans certains cas, l'intégralité des coefficients de Taylor de racines d'applications miroirs et de type miroir. En particulier, nous démontrons une conjecture faite par Zhou dans \cite{Zhou1} et nous vérifions en partie une observation faite par Krattenthaler et Rivoal dans \cite{Tanguy1}.

Dans un premier temps, nous énonçons les deux résultats principaux de \cite{Delaygue}. Pour cela, nous devons introduire quelques notations.

Dans la suite, si $\textbf{e}:=(e_1,\dots,e_{q_1})$ et $\textbf{f}:=(f_1,\dots,f_{q_2})$ sont deux suites d'entiers positifs, on note $|\textbf{e}|:=\sum_{i=1}^{q_1}e_i$, $|\textbf{f}\,|:=\sum_{j=1}^{q_2}f_j$, $M_{\textbf{e},\textbf{f}}:=\max(e_1,\dots,e_{q_1},f_1,\dots,f_{q_2})$ et 
$$
\mathcal{Q}_{\textbf{e},\textbf{f}}(n):=\frac{(e_1n)!\cdots(e_{q_1}n)!}{(f_1n)!\cdots(f_{q_2}n)!},
$$ 
où $n\in\mathbb{N}$. On définit les séries entières
$$
F_{\textbf{e},\textbf{f}}(z):=\sum_{n=0}^{\infty}\frac{(e_1n)!\cdots(e_{q_1}n)!}{(f_1n)!\cdots(f_{q_2}n)!}z^n,
$$
\begin{equation}\label{definition de G facto}
G_{\textbf{e},\textbf{f}}(z):=\sum_{n=1}^{\infty}\frac{(e_1n)!\cdots(e_{q_1}n)!}{(f_1n)!\cdots(f_{q_2}n)!}\left(\sum_{i=1}^{q_1}e_iH_{e_in}-\sum_{j=1}^{q_2}f_jH_{f_jn}\right)z^n
\end{equation}
et, pour tout $L\in\{1,\dots,M_{\textbf{e},\textbf{f}}\}$, on pose
\begin{equation}\label{définition G L}
G_{L,\textbf{e},\textbf{f}}(z):=\sum_{n=1}^{\infty}\frac{(e_1n)!\cdots(e_{q_1}n)!}{(f_1n)!\cdots(f_{q_2}n)!}H_{Ln}z^n,
\end{equation}
où $H_n:=\sum_{i=1}^n\frac{1}{i}$ est le $n$-ième nombre harmonique. La série $F_{\textbf{e},\textbf{f}}(z)$ est une série hypergéométrique généralisée et est donc solution d'une équation différentielle linéaire. On obtient, \textit{via} la méthode de Frobenius (voir \cite{Schwarz}), une base de solutions de cette équation différentielle avec, dans certains cas (voir \cite[Partie $7$]{Delaygue}), au plus des singularités logarithmiques à l'origine, dont $F_{\textbf{e},\textbf{f}}(z)$ et $G_{\textbf{e},\textbf{f}}(z)+\log(z)F_{\textbf{e},\textbf{f}}(z)$.
 
Dans le contexte de la symétrie miroir, cas où $|\textbf{e}|=|\textbf{f}\,|$, la fonction 
$$
q_{\textbf{e},\textbf{f}}(z):=z\exp(G_{\textbf{e},\textbf{f}}(z)/F_{\textbf{e},\textbf{f}}(z))
$$ 
est une \textit{coordonnée canonique} et son inverse pour la composition $z(q)$ est une \textit{application miroir}(\footnote{Dans cet article, nous énoncerons uniquement des résultats d'intégralité pour les coefficients de Taylor de racines de coordonnées canoniques. Ils vaudront tout autant pour les racines d'applications miroir puisque, d'après l'introduction de \cite{Lian Yau 2}, pour tout $v\in\mathbb{N}$, $v\geq 1$, on a $(q^{-1}z(q))^{1/v}\in \mathbb{Z}[[q]]$ si, et seulement si $(z^{-1}q(z))^{1/v}\in\mathbb{Z}[[z]]$.}). Pour tout $L\in\{1,\dots,M_{\textbf{e},\textbf{f}}\}$, on définit $q_{L,\textbf{e},\textbf{f}}(z):=\exp(G_{L,\textbf{e},\textbf{f}}(z)/F_{\textbf{e},\textbf{f}}(z))$. On dira que $q_{L,\textbf{e},\textbf{f}}$ est une application de \textit{type miroir}. On a alors la relation 
$$
q_{\textbf{e},\textbf{f}}(z):=z\frac{\prod_{i=1}^{q_1}q_{e_i,\textbf{e},\textbf{f}}^{e_i}(z)}{\prod_{i=1}^{q_2}q_{f_i,\textbf{e},\textbf{f}}^{f_i}(z)}.
$$

Afin d'énoncer les critères démontrés dans \cite{Delaygue}, nous devons introduire l'\textit{application de Landau} $\Delta_{\textbf{e},\textbf{f}}$ définie comme suit : pour tout $x\in\mathbb{R}$, on pose 
$$
\Delta_{\textbf{e},\textbf{f}}(x):=\sum_{i=1}^{q_1}\lfloor e_ix\rfloor-\sum_{i=1}^{q_2}\lfloor f_ix\rfloor,
$$
où $\lfloor\cdot\rfloor$ désigne la fonction partie entière. On note également $\{\cdot\}$ la fonction partie fractionnaire. La fonction $\Delta_{\textbf{e},\textbf{f}}$ est constante par morceaux. Pour tout entier positif $c$, on a $\lfloor cx\rfloor=\lfloor c\{x\}\rfloor+c\lfloor x\rfloor$ et donc $\Delta_{\textbf{e},\textbf{f}}(x)=\Delta_{\textbf{e},\textbf{f}}(\{x\})+(|\textbf{e}|-|\textbf{f}\,|)\lfloor x\rfloor$. Ainsi, on a $|\textbf{e}|=|\textbf{f}\,|$ si et seulement si $\Delta_{\textbf{e},\textbf{f}}$ est $1$-périodique. Si $M_{\textbf{e},\textbf{f}}=0$, alors $\Delta_{\textbf{e},\textbf{f}}$ est la fonction nulle. Sinon, la fonction $\Delta_{\textbf{e},\textbf{f}}$ est nulle sur $[0,1/M_{\textbf{e},\textbf{f}}[$ car $0\leq e_i/M_{\textbf{e},\textbf{f}}\leq 1$ et $0\leq f_i/M_{\textbf{e},\textbf{f}}\leq 1$. 

On rappelle le \textit{Critère de Landau} (voir \cite{Landau}) (\footnote{L'énoncé présenté dans cet article est le critère de Landau de \cite{Landau} appliqué en une dimension.}) : si $\textbf{e}$ et $\textbf{f}$ sont deux suites d'entiers positifs, alors tous les termes de la suite $\mathcal{Q}_{\textbf{e},\textbf{f}}$ sont entiers si, et seulement si, pour tout $x\in[0,1]$, on a $\Delta_{\textbf{e},\textbf{f}}(x)\geq 0$. Nous pouvons maintenant énoncer les deux résultats principaux de \cite{Delaygue}.

\begin{Crit1}[Théorème $1$ de \cite{Delaygue}]\label{conj equiv}
Soit $\textbf{e}$ et $\textbf{f}$ deux suites d'entiers strictement positifs disjointes telles que $\mathcal{Q}_{\textbf{e},\textbf{f}}$ soit une suite à termes entiers (ce qui équivaut à $\Delta_{\textbf{e},\textbf{f}}\geq 0$ sur $[0,1]$) et vérifiant $|\textbf{e}|=|\textbf{f}\,|$. On a alors la dichotomie suivante.
\begin{itemize}
\item[$(i)$] Si, pour tout $x\in[1/M_{\textbf{e},\textbf{f}},1[$, on a $\Delta_{\textbf{e},\textbf{f}}(x)\geq 1$, alors $q_{\textbf{e},\textbf{f}}(z)\in z\mathbb{Z}[[z]]$.
\item[$(ii)$] S'il existe $x\in[1/M_{\textbf{e},\textbf{f}},1[$ tel que $\Delta_{\textbf{e},\textbf{f}}(x)=0$, alors il n'existe qu'un nombre fini de nombres premiers~$p$ tels que $q_{\textbf{e},\textbf{f}}(z)\in z\mathbb{Z}_p[[z]]$.
\end{itemize}
\end{Crit1}

Les applications de type miroir ont initialement été introduites par Krattenthaler et Rivoal dans \cite{Tanguy1} pour simplifier la preuve de l'intégralité des coefficients de Taylor des applications miroir lorsque $\Delta_{\textbf{e},\textbf{f}}$ est croissante sur $[0,1[$ et $1$-périodique. On dispose aussi d'un critère d'intégralité pour leurs coefficients de Taylor.

\begin{Crit2}[Théorème $2$ de \cite{Delaygue}]\label{conj equiv L}
Soit $\textbf{e}$ et $\textbf{f}$ deux suites d'entiers strictement positifs disjointes telles que $\mathcal{Q}_{\textbf{e},\textbf{f}}$ soit une suite à termes entiers (ce qui équivaut à $\Delta_{\textbf{e},\textbf{f}}\geq 0$ sur $[0,1]$) et vérifiant $|\textbf{e}|=|\textbf{f}\,|$. On a alors la dichotomie suivante.
\begin{itemize}
\item[$(i)$] Si, pour tout $x\in[1/M_{\textbf{e},\textbf{f}},1[$, on a $\Delta_{\textbf{e},\textbf{f}}(x)\geq 1$, alors, pour tout $L\in\{1,\dots,M_{\textbf{e},\textbf{f}}\}$, on a $q_{L,\textbf{e},\textbf{f}}(z)\in \mathbb{Z}[[z]]$.
\item[$(ii)$] S'il existe $x\in[1/M_{\textbf{e},\textbf{f}},1[$ tel que $\Delta_{\textbf{e},\textbf{f}}(x)=0$, alors, pour tout $L\in\{1,\dots,M_{\textbf{e},\textbf{f}}\}$, il n'existe qu'un nombre fini de nombres premiers~$p$ tels que $q_{L,\textbf{e},\textbf{f}}(z)\in \mathbb{Z}_p[[z]]$.
\end{itemize}
\end{Crit2}

Pour espérer qu'une racine de $z^{-1}q_{\textbf{e},\textbf{f}}(z)$ ou qu'une racine d'une des $q_{L,\textbf{e},\textbf{f}}(z)$ ait tous ses coefficients de Taylor entiers, ces deux critères montrent qu'il faut que, pour tout $x\in[1/M_{\textbf{e},\textbf{f}},1[$, on ait $\Delta_{\textbf{e},\textbf{f}}(x)\geq 1$. Avant d'énoncer nos résultats sur l'intégralité des coefficients de Taylor de racines d'applications miroir ou de type miroir, nous donnons une précision sur les éventuels $v\in\mathbb{N}$, $v\geq 1$, vérifiant $(z^{-1}q_{\textbf{e},\textbf{f}}(z))^{1/v}\in\mathbb{Z}[[z]]$ ou $q_{L,\textbf{e},\textbf{f}}(z)^{1/v}\in\mathbb{Z}[[z]]$. Soit $q(z)\in\mathbb{Z}[[z]]$ et $V$ le plus grand entier naturel tel que $q(z)^{1/V}\in\mathbb{Z}[[z]]$. Le lemme $5$ de \cite{Heninger} donne que les entiers naturels $U$ tels que $q(z)^{1/U}\in\mathbb{Z}[[z]]$ sont exactement les diviseurs positifs de $V$.

Dans la suite de l'article, si $\textbf{e}$ et $\textbf{f}$ sont deux suites d'entiers positifs, alors, pour tout élément $L$ de $\textbf{e}$ et $\textbf{f}$, on note $D_L$ le plus petit multiple commun des entiers allant de $1$ à $\lfloor M_{\textbf{e},\textbf{f}}/L\rfloor$. Le théorème suivant raffine le point $(i)$ du critère pour les applications de type miroir.

\begin{theo}\label{RootsL}
Soit $\textbf{e}$ et $\textbf{f}$ deux suites d'entiers strictement positifs disjointes vérifiant $|\textbf{e}|=|\textbf{f}\,|$ et telles que, pour tout $x\in[1/M_{\textbf{e},\textbf{f}},1[$, on a $\Delta_{\textbf{e},\textbf{f}}(x)\geq 1$. Alors, pour tout $L\in\{1,\dots,M_{\textbf{e},\textbf{f}}\}$, on a 
$$
q_{L,\textbf{e},\textbf{f}}(z)^{1/D_L}\in\mathbb{Z}[[z]].
$$
\end{theo}

\begin{Remarque}
Le théorème \ref{RootsL} est principalement dû à une généralisation du lemme $9$ de \cite{Delaygue} qui est ici le lemme \ref{lemme 24 généralisé} de la partie \ref{demo lemme 10}.
\end{Remarque}

Le théorème \ref{RootsL} permet, dans certains cas, d'obtenir l'intégralité des coefficients de Taylor de racines d'applications miroir. Dans la suite, pour tout $a$ et $b$ entiers, on note $\textup{pgcd}(a,b)$ le plus grand diviseur commun de $a$ et $b$.

\begin{theo}\label{cor1}
Soit $\textbf{e}$ et $\textbf{f}$ deux suites d'entiers strictement positifs disjointes vérifiant $|\textbf{e}|=|\textbf{f}\,|$ et telles que, pour tout $x\in[1/M_{\textbf{e},\textbf{f}},1[$, on a $\Delta_{\textbf{e},\textbf{f}}(x)\geq 1$. Soit $\theta\geq 1$ un diviseur de $M_{\textbf{e},\textbf{f}}$. Si, pour tout élément $L$ de $\textbf{e}$ et $\textbf{f}$, l'entier $\theta/\textup{pgcd}(L,\theta)$ divise $D_L$, alors on a
$$
(z^{-1}q_{\textbf{e},\textbf{f}}(z))^{1/\theta}\in\mathbb{Z}[[z]].
$$ 
\end{theo}

Dans la partie \ref{section énoncé}, nous montrerons que les théorèmes \ref{RootsL} et \ref{cor1} ne sont pas optimaux dans le sens où il existe des suites $\textbf{e}$ et $\textbf{f}$ pour lesquelles ces théorèmes ne donnent pas les entiers maximaux $V$ et $V_L$ tels que $(z^{-1}q_{\textbf{e},\textbf{f}}(z))^{1/V}\in\mathbb{Z}[[z]]$ et $q_{L,\textbf{e},\textbf{f}}(z)^{1/V_L}\in\mathbb{Z}[[z]]$. 

Nous donnons maintenant une certaine classe de suites $\textbf{e}$ et $\textbf{f}$ pour lesquelles on peut appliquer le théorème \ref{cor1} avec $M_{\textbf{e},\textbf{f}}$ à la place de $\theta$. En particulier, nous démontrons la conjecture de Zhou \cite[Conjecture 1]{Zhou1} qui s'énonce comme suit. 

Soient $k_1,\dots,k_n$ des entiers strictement positifs vérifiant $\frac{1}{k_1}+\cdots+\frac{1}{k_n}=1$. Soit $k$ le plus petit multiple commun de $k_1,\dots,k_n$ et, pour tout $i\in\{1,\dots,n\}$, $w_i:=k/k_i$. Alors, en notant $\textbf{e}:=(k)$ et $\textbf{f}:=(w_1,\dots,w_n)$, on a 
\begin{equation}\label{corZhou}
(z^{-1}q_{\textbf{e},\textbf{f}}(z))^{1/k}\in\mathbb{Z}[[z]].
\end{equation}

\begin{cor}\label{cor2}
Soit $\textbf{e}$ et $\textbf{f}$ deux suites d'entiers strictement positifs disjointes vérifiant $|\textbf{e}|=|\textbf{f}\,|$ telles que, pour tout $x\in[1/M_{\textbf{e},\textbf{f}},1[$, on a $\Delta_{\textbf{e},\textbf{f}}(x)\geq 1$ et telles que tout élément de $\textbf{e}$ et $\textbf{f}$ divise $M_{\textbf{e},\textbf{f}}$. Alors, on a
$$
(z^{-1}q_{\textbf{e},\textbf{f}}(z))^{1/M_{\textbf{e},\textbf{f}}}\in\mathbb{Z}[[z]].
$$
En particulier, la conjecture de Zhou est vraie.
\end{cor}

Dans la partie \ref{section énoncé}, nous remarquons que les résultats antérieurs connus sur l'intégralité des coefficients de Taylor de racines d'applications miroir ou de type miroir s'appliquent uniquement à des suites $\textbf{e}$ et $\textbf{f}$ telles que $\Delta_{\textbf{e},\textbf{f}}$ est croissante sur $[0,1[$ et $1$-périodique.
Nous montrons sur un exemple que les suites $\textbf{e}$ et $\textbf{f}$ traitées par la conjecture de Zhou ne définissent pas forcément une application $\Delta_{\textbf{e},\textbf{f}}$ croissante sur $[0,1[$. 

En effet, prenons $k_1=3$, $k_2=k_3=4$ et $k_4=6$, on obtient $\frac{1}{k_1}+\frac{1}{k_2}+\frac{1}{k_3}+\frac{1}{k_4}=\frac{1}{3}+\frac{1}{2}+\frac{1}{6}=1$. On a $k=12$, $w_1=4$, $w_2=w_3=3$ et $w_4=2$. On utilise maintenant la méthode expliquée dans la partie $7$ de \cite{Delaygue} pour obtenir les sauts de $\Delta_{\textbf{e},\textbf{f}}$ sur $[0,1]$. Pour cela, pour tout $x\in\mathbb{R}$ et tout $n\in\mathbb{N}$, on note $(x)_n:=x(x+1)\cdots(x+n-1)$ pour $n\geq 1$ et $(x)_0=1$ (symbole de Pochhammer). On écrit
\begin{align*}
\frac{(12n)!}{(4n)!(3n)!^2(2n)!}
&=\left(\frac{12^{12}}{4^4(3^3)^22^2}\right)^n\frac{(1/12)_n(2/12)_n\cdots(12/12)_n}{(1/4)_n(2/4)_n(3/4)_n(4/4)_n(1/3)_n^2(2/3)_n^2(3/3)_n^2(1/2)_n(2/2)_n}\\
&=\left(\frac{12^{12}}{4^43^62^2}\right)^n\frac{(1/12)_n(1/6)_n(5/12)_n(7/12)_n(5/6)_n(11/12)_n}{(1/3)_n(1/2)_n(2/3)_n(1)_n^2}
\end{align*}
et on en déduit que $\Delta_{\textbf{e},\textbf{f}}$ a des sauts d'amplitude $-1$ aux abscisses $1/3$, $1/2$ et $2/3$ donc $\Delta_{\textbf{e},\textbf{f}}$ n'est pas croissante sur $[0,1[$.

\subsection{\'Enoncés des résultats antérieurs et comparaison avec les théorèmes \ref{RootsL} et~\ref{cor1}}\label{section énoncé}

Le premier résultat sur l'intégralité des coefficients de Taylor de racines d'applications miroir est dû à Lian et Yau et s'énonce comme suit.

\begin{theo2}[Lian, Yau, \cite{Lian Yau 1}]
Soit $p\geq 3$ un nombre premier, $\textbf{e}=(p)$ et $\textbf{f}=(1,\dots,1)$ avec $|\textbf{e}|=|\textbf{f}\,|$. Alors, on a 
$$
(z^{-1}q_{\textbf{e},\textbf{f}}(z))^{1/p}\in\mathbb{Z}[[z]].
$$
\end{theo2}

Ce théorème est contenu dans le corollaire \ref{cor2}. Krattenthaler et Rivoal se sont intéressés, dans \cite{Tanguy1} puis dans \cite{Tanguy2}, à l'intégralité des coefficients de Taylor de racines d'applications miroir et de type miroir associées à des suites $\textbf{e}$ et $\textbf{f}$ d'une forme particulière. Nous énonçons leurs résultats par généralité décroissante relativement à la forme des suites $\textbf{e}$ et $\textbf{f}$. Le premier résultat (cf. \cite[Remarque $(b)$, p.181]{Tanguy1}) peut se reformuler comme suit (\footnote{Nous expliquons, à la fin de la partie $7.2$ de \cite{Delaygue} que les suites $\textbf{e}$ et $\textbf{f}$ considérées dans la remarque $(b)$ de la page $181$ de \cite{Tanguy1} sont exactement celles pour lesquelles l'application $\Delta_{\textbf{e},\textbf{f}}$ est croissante sur $[0,1[$ et $1$-périodique.}).

\begin{theo2}[Krattenthaler, Rivoal, \cite{Tanguy1}]
Soit $\textbf{e}$ et $\textbf{f}$ deux suites d'entiers positifs vérifiant $|\textbf{e}|=|\textbf{f}\,|$. Si $\Delta_{\textbf{e},\textbf{f}}$ est croissante sur $[0,1[$, alors le plus grand entier naturel $V$ tel que $q_{1,\textbf{e},\textbf{f}}(z)^{1/V}\in\mathbb{Z}[[z]]$ est $\mathcal{Q}_{\textbf{e},\textbf{f}}(1)$.
\end{theo2}

Par exemple, en prenant $\textbf{e}=(6)$ et $\textbf{f}=(3,2,1)$, on obtient que $\Delta_{\textbf{e},\textbf{f}}$ est croissante sur $[0,1[$ et que $\mathcal{Q}_{\textbf{e},\textbf{f}}(1)=\frac{6!}{3!2!}=60$. D'après le théorème précédent, on obtient que $q_{1,\textbf{e},\textbf{f}}(z)^{1/60}\in\mathbb{Z}[[z]]$. Krattenthaler et Rivoal remarquent dans \cite{Tanguy1} qu'il semblerait que les relations suivantes soient les meilleures possibles : $q_{2,\textbf{e},\textbf{f}}(z)^{1/6}$, $q_{3,\textbf{e},\textbf{f}}(z)^{1/2}$, $q_{4,\textbf{e},\textbf{f}}(z)$, $q_{5,\textbf{e},\textbf{f}}(z)$ et $q_{6,\textbf{e},\textbf{f}}(z)$ sont dans $\mathbb{Z}[[z]]$. En appliquant le théorème \ref{RootsL} avec $\textbf{e}=(6)$ et $\textbf{f}=(3,2,1)$, on obtient bien que $q_{1,\textbf{e},\textbf{f}}(z)^{1/60}$, $q_{2,\textbf{e},\textbf{f}}(z)^{1/6}$, $q_{3,\textbf{e},\textbf{f}}(z)^{1/2}$, $q_{4,\textbf{e},\textbf{f}}(z)$, $q_{5,\textbf{e},\textbf{f}}(z)$ et $q_{6,\textbf{e},\textbf{f}}(z)$ sont dans $\mathbb{Z}[[z]]$. En appliquant le corollaire \ref{cor2} avec $\textbf{e}=(6)$ et $\textbf{f}=(3,2,1)$, on obtient que $(z^{-1}q_{\textbf{e},\textbf{f}}(z))^{1/6}\in\mathbb{Z}[[z]]$.

Dans le cas particulier où $f_1=\cdots=f_{q_2}=1$ (\footnote{D'après la partie $7$ de \cite{Delaygue}, le fait que $f_1=\cdots=f_{q_2}=1$ impose que $\Delta_{\textbf{e},\textbf{f}}$ est croissante sur $[0,1[$.}), Krattenthaler et Rivoal montrent l'intégralité des coefficients de Taylor de racines d'applications de type miroir $q_{L,\textbf{e},\textbf{f}}$ pour tous les $L$ dans $\{1,\dots,M_{\textbf{e},\textbf{f}}\}$, l'exposant donné dans le cas où $L=1$ est l'exposant maximal $\mathcal{Q}_{\textbf{e},\textbf{f}}(1)$.

\begin{theo2}[Krattenthaler, Rivoal, \cite{Tanguy1}]
Soit $\textbf{e}$ et $\textbf{f}:=(1,\dots,1)$ deux suites d'entiers positifs vérifiant $|\textbf{e}|=|\textbf{f}\,|$. Soit $\Theta_L:=L!/\textup{pgcd}(L!,L!H_L)$ le dénominateur de $H_L$ écrit comme une fraction irréductible. Alors, pour tout $L\in\{1,\dots,M_{\textbf{e},\textbf{f}}\}$, on a 
$$
q_{L,\textbf{e},\textbf{f}}(z)^{\Theta_L/\mathcal{Q}_{\textbf{e},\textbf{f}}(1)}\in\mathbb{Z}[[z]].
$$
\end{theo2}

On remarque que, pour tout $L\in\{1,\dots,M_{\textbf{e},\textbf{f}}\}$, on a $\mathcal{Q}_{\textbf{e},\textbf{f}}(1)/\Theta_L\in\mathbb{N}$. En effet, pour tout $L\in\{1,\dots,M_{\textbf{e},\textbf{f}}\}$, $\Theta_L$ divise $L!$ qui divise $e_1!\cdots e_{q_1}!=\mathcal{Q}_{\textbf{e},\textbf{f}}(1)$. 

L'exposant donné par ce théorème est meilleur que celui que l'on obtient avec le théorème~\ref{RootsL} car, pour tout $L\in\{1,\dots,M_{\textbf{e},\textbf{f}}\}$, $D_L$ divise $\mathcal{Q}_{\textbf{e},\textbf{f}}(1)/\Theta_L$. En effet, soit $p$ un nombre premier et $\alpha:=v_p(D_L)$. On a alors $p^{\alpha}\leq\lfloor M_{\textbf{e},\textbf{f}}/L\rfloor$, \textit{i.e.} $M_{\textbf{e},\textbf{f}}\geq Lp^{\alpha}$ et on obtient bien que $p^{\alpha}$ divise $M_{\textbf{e},\textbf{f}}!/L!$ qui divise $M_{\textbf{e},\textbf{f}}!/\Theta_L$ qui enfin divise $\mathcal{Q}_{\textbf{e},\textbf{f}}(1)/\Theta_L$.

Dans le cas où $e_1=\cdots=e_{q_1}=:N\geq 2$ et $f_1=\cdots=f_{q_2}=1$, Krattenthaler et Rivoal améliorent l'exposant obtenu \textit{via} le théorème précédent pour $L=N$, c'est le théorème $1$ de \cite{Tanguy2} qui s'énonce comme suit.

\begin{theo2}[Krattenthaler, Rivoal, \cite{Tanguy2}]
Soit $N\geq 2$ un entier. Soit $\textbf{e}:=(N,\dots,N)$ et $\textbf{f}:=(1,\dots,1)$ vérifiant $|\textbf{e}|=|\textbf{f}\,|$. Soit
$$
\Xi_N:=\prod_{p\leq N}p^{\min\{2+\xi(p,N),v_p(H_N)\}},
$$
où $\xi(p,N)=1$ si $p$ est un nombre premier de Wolstenholme (\footnote{Un nombre premier de Wolstenholme est un premier $p$ qui vérifie $v_p(H_{p-1})\geq 3$. Actuellement, les seuls premiers de Wolstenholme connus sont $16843$ et $2124679$, et on ne sait pas si l'ensemble des premiers de Wolstenholme est fini ou infini.}) ou si $p$ divise $N$, et $\xi(p,N)=0$ sinon. Alors, on a
$$
q_{N,\textbf{e},\textbf{f}}(z)^{\frac{1}{\Xi_N\mathcal{Q}_{\textbf{e},\textbf{f}}(1)}}\in\mathbb{Z}[[z]].
$$
\end{theo2}

On remarque que $\Xi_N$ n'est pas forcément un entier mais que $\Xi_N\mathcal{Q}_{\textbf{e},\textbf{f}}(1)$ en est un. 
En effet, si on regroupe, dans la définition de $\Xi_N$, les premiers $p$ pour lesquels on a $v_p(H_N)\leq -1$, on obtient que $\Xi_N=w_N/\Theta_N$ où $w_N\in\mathbb{N}$. Ainsi, on a $\Xi_N\mathcal{Q}_{\textbf{e},\textbf{f}}(1)=\mathcal{Q}_{\textbf{e},\textbf{f}}(1)w_N/\Theta_N$ avec, comme on l'a vu précédemment, $\mathcal{Q}_{\textbf{e},\textbf{f}}(1)/\Theta_N\in\mathbb{N}$. En particulier, on voit que ce théorème améliore le théorème précédent d'un facteur $w_N$ lorsque $\textbf{e}=(N,\dots,N)$ et $L=N$.

Il est conjecturé dans \cite{Tanguy2} que ce théorème est optimal lorsque $\textbf{e}=(N)$, \textit{i.e.} que dans ce cas, le plus grand entier $V_N$ tel que $q_{N,\textbf{e},\textbf{f}}(z)^{1/V_N}\in\mathbb{Z}[[z]]$ est $V_N=\Xi_N\mathcal{Q}_{\textbf{e},\textbf{f}}(1)=\Xi_NN!$. Il est expliqué, toujours dans \cite{Tanguy2}, que cette conjecture est impliquée par la conjecture sur les nombres harmoniques suivante : il n'existe pas de premier $p$ et d'entier $N\geq 1$ tels que $v_p(H_N)\geq 4$.

Toujours dans le cas où $e_1=\cdots=e_{q_1}=N\geq 2$ et $f_1=\cdots=f_{q_2}=1$, le théorème $2$ de \cite{Tanguy2} donne l'intégralité des coefficients de Taylor de racines d'applications miroir :

\begin{theo2}[Krattenthaler, Rivoal, \cite{Tanguy2}]
Soit $N\geq 2$ un entier. Soit $\textbf{e}=(N,\dots,N)$ avec $q_1$ occurences de $N$, et $\textbf{f}=(1,\dots,1)$ vérifiant $|\textbf{e}|=|\textbf{f}\,|$. Soit
$$
\Omega_N:=\prod_{p\leq N}p^{\min\{2+\omega(p,N),v_p(H_N-1)\}},
$$
où $\omega(p,N)=1$ si $p$ est un nombre premier de Wolstenholme ou si $N\equiv\pm 1\mod p$, et $\omega(p,N)=0$ sinon. Alors, on a
$$
(z^{-1}q_{\textbf{e},\textbf{f}}(z))^{\frac{1}{\Omega_N\mathcal{Q}_{\textbf{e},\textbf{f}}(1)q_1N}}\in\mathbb{Z}[[z]].
$$
\end{theo2}

On remarque que $\Omega_N$ n'est pas forcément un entier mais que $\Omega_N\mathcal{Q}_{\textbf{e},\textbf{f}}(1)$ en est un. En effet, si on regroupe, dans la définition de $\Omega_N$, les premiers $p$ pour lesquels on a $v_p(H_N-1)\leq -1$, on obtient que $\Omega_N=m_N/\Theta_N$, où $m_N\in\mathbb{N}$, car $H_N$ et $H_N-1$ ont le même dénominateur, une fois écrits sous formes irréductibles. Ainsi, on a $\Omega_N\mathcal{Q}_{\textbf{e},\textbf{f}}(1)=\mathcal{Q}_{\textbf{e},\textbf{f}}(1)m_N/\Theta_N$ avec, comme on l'a vu précédemment, $\mathcal{Q}_{\textbf{e},\textbf{f}}(1)/\Theta_N\in\mathbb{N}$. En particulier, on voit que dans le cas où $\textbf{e}=(N,\dots,N)$ et $\textbf{f}=(1,\dots,1)$, l'exposant donné par ce théorème est nettement meilleur que l'exposant $\theta$ obtenu par le théorème \ref{cor1} puisque $\theta$ est un diviseur de $M_{\textbf{e},\textbf{f}}=N$.

Il est conjecturé dans \cite{Tanguy2} que le théorème ci-dessus est optimal lorsque $\textbf{e}=(N)$, \textit{i.e.} que dans ce cas, le plus grand entier $V_N$ tel que $(z^{-1}q_{\textbf{e},\textbf{f}}(z))^{1/V_N}\in\mathbb{Z}[[z]]$ est $V_N=\Omega_N\mathcal{Q}_{\textbf{e},\textbf{f}}(1)N=\Omega_NN!N$. Il est expliqué, toujours dans \cite{Tanguy2}, que cette conjecture est impliquée par la conjecture sur les nombres harmoniques suivante : il n'existe pas de premier $p$ et d'entier $N\geq 2$ tels que $v_p(H_N-1)\geq 4$.

\subsection{Trame de la démonstration des théorèmes \ref{RootsL} et \ref{cor1}}

Dans la partie \ref{section equiv Zp}, on ramène le théorème \ref{RootsL} à la preuve d'un énoncé $p$-adique. 

La partie \ref{section sur la condition suffisante} est consacrée à la suite de la preuve du théorème \ref{RootsL} qui contient la principale amélioration de la preuve du point $(i)$ du critère pour les applications de type miroir.

Dans la partie \ref{partie4}, on démontre le théorème \ref{cor1} et le corollaire \ref{cor2} en utilisant le théorème \ref{RootsL}.

\section{Une réduction $p$-adique du théorème \ref{RootsL}}\label{section equiv Zp}

On se place sous les hypothèses du théorème \ref{RootsL} et on fixe $L\in\{1,\dots,M_{\textbf{e},\textbf{f}}\}$ dans cette partie. Pour alléger les notations, on notera $\Delta:=\Delta_{\textbf{e},\textbf{f}}$, $\mathcal{Q}:=\mathcal{Q}_{\textbf{e},\textbf{f}}$, $F:=F_{\textbf{e},\textbf{f}}$, $G:=G_{\textbf{e},\textbf{f}}$, $G_L:=G_{L,\textbf{e},\textbf{f}}$, $q:=q_{\textbf{e},\textbf{f}}$ et $q_L:=q_{L,\textbf{e},\textbf{f}}$, comme dans toute la suite de cet article. 

On rappelle que $q_{L}(z)^{1/D_L}\in\mathbb{Z}[[z]]$, si, et seulement si, pour tout nombre premier $p$, on a $q_{L}(z)^{1/D_L}\in\mathbb{Z}_p[[z]]$. 

Nous allons définir, pour tout nombre premier $p$, des éléments $\Phi_{L,p}(a+Kp)$ de $\mathbb{Q}_p$, où $a\in\{0,\dots,p-1\}$ et $K\in\mathbb{N}$, et montrer que $q_{L}(z)^{1/D_L}\in\mathbb{Z}[[z]]$, si, et seulement si, pour tout nombre premier $p$, tout $a\in\{0,\dots,p-1\}$ et tout $K\in\mathbb{N}$, on a $\Phi_{L,p}(a+Kp)\in pD_L\mathbb{Z}_p$. 

On fixe un nombre premier $p$ dans cette partie.

Avant de donner la preuve du théorème \ref{RootsL}, on va le reformuler. Le résultat classique suivant est dû à Dieudonné et Dwork (voir \cite[Chap. VI, Sec. 2, Lemma 3]{Kobliz}; \cite[Chap. 14, Sec. 2]{Lang}).

\begin{lemme}\label{réduction modulo p}
Soit $F(z)$ une série formelle dans $1+z\mathbb{Q}[[z]]$. Alors $F(z)\in1+z\mathbb{Z}_p[[z]]$ si et seulement si $\frac{F(z^p)}{F(z)^p}\in 1+pz\mathbb{Z}_p[[z]]$.
\end{lemme}

On déduit de ce lemme le corollaire suivant (voir \cite[Lemma 5, p. 610]{Zudilin}), qui va nous permettre \og d'éliminer\fg\,l'exponentielle dans l'expression $q_L(z)^{1/D_L}=\exp(G_L(z)/D_LF(z))$.

\begin{cor}\label{cor reduction}
Soit $f(z)\in z\mathbb{Q}[[z]]$. On a $e^{f(z)}\in 1+z\mathbb{Z}_p[[z]]$ si et seulement si $f(z^p)-pf(z)\in pz\mathbb{Z}_p[[z]]$.
\end{cor}

D'après l'identité \eqref{définition G L} définissant $G_L$, on a $G_L(0)=0$ et donc $G_L(z)/D_LF(z)$ est dans $z\mathbb{Q}[[z]]$. Ainsi, d'après le corollaire \ref{cor reduction}, on a $q_L(z)^{1/D_L}\in\mathbb{Z}_p[[z]]$, si, et seulement si $\frac{G_L}{F}(z^p)-p\frac{G_L}{F}(z)\in pD_Lz\mathbb{Z}_p[[z]]$. 

Or, comme $\mathcal{Q}$ est une suite à termes entiers, on a $F(z)\in 1+z\mathbb{Z}[[z]]\subset 1+z\mathbb{Z}_p[[z]]$. Ainsi, $q_L(z)^{1/D_L}\in\mathbb{Z}_p[[z]]$, si, et seulement si on a $F(z)G_L(z^p)-pF(z^p)G_L(z)\in pD_Lz\mathbb{Z}_p[[z]]$.

D'après l'identité \eqref{définition G L} définissant $G_L$, le coefficient de $z^{a+Kp}$ dans $F(z)G_L(z^p)-pF(z^p)G_L(z)$ est 
$$
\Phi_{L,p}(a+Kp):=\sum_{j=0}^K\mathcal{Q}(K-j)\mathcal{Q}(a+jp)(H_{L(K-j)}-pH_{L(a+jp)}).
$$

On a donc $q_L(z)^{1/D_L}\in\mathbb{Z}_p[[z]]$, si, et seulement si, pour tout $a\in\{0,\dots,p-1\}$ et tout $K\in\mathbb{N}$, on a $\Phi_{L,p}(a+Kp)\in pD_L\mathbb{Z}_p$.

\section{Démonstration du théorème \ref{RootsL}}\label{section sur la condition suffisante}

On se place sous les hypothèses du thèorème \ref{RootsL}. Le but de cette partie est de montrer que, pour tout $L\in\{1,\dots,M_{\textbf{e},\textbf{f}}\}$, on a $q_L(z)^{1/D_L}\in\mathbb{Z}[[z]]$. D'après la partie \ref{section equiv Zp}, il nous suffit de montrer que, pour tout $L\in\{1,\dots,M_{\textbf{e},\textbf{f}}\}$, tout nombre premier $p$, tout $a\in\{0,\dots,p-1\}$ et tout $K\in\mathbb{N}$, on a $\Phi_{L,p}(a+Kp)\in pD_L\mathbb{Z}_p$. On fixe $L\in\{1,\dots,M_{\textbf{e},\textbf{f}}\}$ dans cette partie.

\subsection{Nouvelle reformulation du problème}\label{newreforprob}

Pour tout premier $p$, tout $a\in\{0,\dots,p-1\}$ et tout $K,s$ et $m$ dans $\mathbb{N}$, on définit 
$$
S(a,K,s,p,m):=\sum_{j=mp^s}^{(m+1)p^s-1}\left(\mathcal{Q}(a+jp)\mathcal{Q}(K-j)-\mathcal{Q}(j)\mathcal{Q}(a+(K-j)p)\right),
$$
où l'on pose $\mathcal{Q}(\ell)=0$ si $\ell$ est un entier strictement négatif.

Le but de cette partie est de produire, pour tout nombre premier $p$, une fonction $g_p$ de $\mathbb{N}$ dans $\mathbb{Z}_p$ telle que: si pour tout premier $p$, tout $a\in\{0,\dots,p-1\}$ et tout $K,s$ et $m$ dans $\mathbb{N}$, on a $S(a,K,s,p,m)\in p^{s+1}g_p(m)\mathbb{Z}_p$, alors on a $\Phi_{L,p}(a+pK)\in pD_L\mathbb{Z}_p$.
Démontrer le théorème~\ref{RootsL} reviendra alors à minorer convenablement la valuation $p$-adique de $S(a,K,s,p,m)$ pour tout nombre premier $p$. Cette méthode de réduction est la même que celle que nous avons utilisé dans \cite{Delaygue} et c'est, à l'origine, une adaptation de l'approche du problème faite par Dwork dans \cite{Dwork 1}. Dans le présent article, nous allons utiliser la même fonction $g_p$ que dans \cite{Delaygue}.

\subsubsection{Une réécriture de $\Phi_{L,p}(a+Kp)$ modulo $p\mathbb{Z}_p$}\label{section une réécriture de}

On fixe un nombre premier $p$. Nous allons montrer que
\begin{equation}\label{*}
\Phi_{L,p}(a+Kp)\equiv -\sum_{j=0}^{K}H_{Lj}\big(\mathcal{Q}(a+jp)\mathcal{Q}(K-j)-\mathcal{Q}(j)\mathcal{Q}(a+(K-j)p)\big)\mod\;pD_L\mathbb{Z}_p.
\end{equation}
Si $a=K=0$, alors c'est évident. On suppose donc, dans cette partie, que $a\neq 0$ ou $K\neq 0$.

Pour tout $a\in\{0,\dots,p-1\}$ et tout $j\in\mathbb{N}$, on a
\begin{align}
pH_{L(a+jp)}
&=p\left(\sum_{i=1}^{Ljp}\frac{1}{i}+\sum_{i=1}^{La}\frac{1}{Ljp+i}\right)\notag\\
&\equiv p\left(\sum_{i=1}^{Lj}\frac{1}{ip}+\sum_{i=1}^{\lfloor La/p\rfloor}\frac{1}{Ljp+ip}\right)\mod p\mathbb{Z}_p\notag\\
&\equiv H_{Lj}+\sum_{i=1}^{\lfloor La/p\rfloor}\frac{1}{Lj+i}\mod p\mathbb{Z}_p\label{(25)}.
\end{align}

Nous avons besoin de deux résultats que l'on démontrera plus loin à l'aide du lemme \ref{lemme 24 généralisé} énoncé dans la partie \ref{demo lemme 10} :

Pour tout $L\in\{1,\dots,M_{\textbf{e},\textbf{f}}\}$, tout $a\in\{0,\dots,p-1\}$, tout $K\in\mathbb{N}$ et tout $j\in\{0,\dots,K\}$ tels que $a\neq 0$ ou $K\neq 0$, on a 
\begin{equation}\label{lemme Q rho 1}
\mathcal{Q}(a+jp)\sum_{i=1}^{\lfloor La/p\rfloor}\frac{1}{Lj+i}\in pD_L\mathbb{Z}_p
\end{equation}
et
\begin{equation}\label{tre21}
\mathcal{Q}(K-j)\mathcal{Q}(a+jp)\in D_1\mathbb{Z}_p\subset D_L\mathbb{Z}_p.
\end{equation}
En appliquant \eqref{lemme Q rho 1} et \eqref{tre21} à \eqref{(25)}, on obtient que \begin{equation}\label{trans93}
\mathcal{Q}(K-j)\mathcal{Q}(a+jp)pH_{L(a+jp)}\equiv\mathcal{Q}(K-j)\mathcal{Q}(a+jp)H_{Lj}\mod pD_L\mathbb{Z}_p.
\end{equation}
Ainsi, on a
\begin{align*}
\Phi_{L,p}(a+Kp)
&=\sum_{j=0}^K\mathcal{Q}(K-j)\mathcal{Q}(a+jp)(H_{L(K-j)}-pH_{L(a+jp)})\\
&\equiv\sum_{j=0}^{K}\mathcal{Q}(K-j)\mathcal{Q}(a+jp)(H_{L(K-j)}-H_{Lj})\mod pD_L\mathbb{Z}_p\\
&\equiv -\sum_{j=0}^{K}H_{Lj}\left(\mathcal{Q}(a+jp)\mathcal{Q}(K-j)-\mathcal{Q}(j)\mathcal{Q}(a+(K-j)p)\right)\mod\;pD_L\mathbb{Z}_p,
\end{align*}
ce qui est bien l'équation \eqref{*} attendue.

On utilise maintenant un lemme combinatoire dû à Dwork (voir \cite[Lemma 4.2, p. 308]{Dwork 1}) qui nous permet d'écrire 
$$
\sum_{j=0}^{K}H_{Lj}\left(\mathcal{Q}(a+jp)\mathcal{Q}(K-j)-\mathcal{Q}(j)\mathcal{Q}(a+(K-j)p)\right)=\sum_{s=0}^r\sum_{m=0}^{p^{r+1-s}-1}W_L(a,K,s,p,m),
$$
où $r$ est tel que $K<p^r$, et 
$$
W_L(a,K,s,p,m):=(H_{Lmp^s}-H_{L\lfloor m/p\rfloor p^{s+1}})S(a,K,s,p,m).
$$ 

Si l'on montre que, pour tout $m$ et $s$ dans $\mathbb{N}$, on a $W_L(a,K,s,p,m)\in pD_L\mathbb{Z}_p$, alors on aura bien $\Phi_{L,p}(a+Kp)\in pD_L\mathbb{Z}_p$, comme voulu.

Pour tout $m\in\mathbb{N}$, on pose $\mu_p(m):=\sum_{\ell=1}^{\infty}\textbf{1}_{[1/M_{\textbf{e},\textbf{f}},1[}(\{m/p^{\ell}\})$, où $\textbf{1}_{[1/M_{\textbf{e},\textbf{f}},1[}$ est la fonction caractéristique de $[1/M_{\textbf{e},\textbf{f}},1[$. Pour tout $m\in\mathbb{N}$, on pose $g_p(m):=p^{\mu_p(m)}$. On utilise maintenant le lemme suivant que l'on démontre dans la partie \ref{demo lemme 10}.

\begin{lemme}\label{valuation diff H}
Pour tout nombre premier $p$, tout $L\in\{1,\dots,M_{\textbf{e},\textbf{f}}\}$ et tout $s$ et $m$ dans $\mathbb{N}$, on~a
$$
p^{s+1}g_p(m)\left(H_{Lmp^s}-H_{L\lfloor m/p\rfloor p^{s+1}}\right)\in pD_L\mathbb{Z}_p.
$$
\end{lemme}

D'après le lemme \ref{valuation diff H}, si on montre que, pour tout $a\in\{0,\dots,p-1\}$, tout $K,s$ et $m$ dans $\mathbb{N}$, on a $S(a,K,s,p,m)\in p^{s+1}g_p(m)\mathbb{Z}_p$, alors, on aura $q_L(z)^{1/D_L}\in\mathbb{Z}_p[[z]]$, ce qui est la reformulation annoncée.

Afin de montrer que,  pour tout $a\in\{0,\dots,p-1\}$, tout $K,s$ et $m$ dans $\mathbb{N}$, on a $S(a,K,s,p,m)\in p^{s+1}g_p(m)\mathbb{Z}_p$, il suffit de suivre à l'identique les parties $4.2$, $4.3$ et $4.4$ de \cite{Delaygue}. La preuve du théorème \ref{RootsL} est donc achevée modulo celles des résultats et lemmes que l'on démontre dans la partie \ref{demo lemme 10}.

\subsubsection{Démonstration de \eqref{lemme Q rho 1}, \eqref{tre21} et du lemme \ref{valuation diff H}}\label{demo lemme 10}

Dans un premier temps, on montre que, pour tout $n\in\mathbb{N}$, on a $v_p\left(\mathcal{Q}(n)\right)=\sum_{\ell=1}^{\infty}\Delta\big(\{n/p^{\ell}\}\big)$. En effet, comme $|\textbf{e}|=|\textbf{f}\,|$, on obtient que $\Delta$ est $1$-périodique. De plus, on rappelle que si $m$ est un entier naturel, on a la formule $v_p(m!)=\sum_{\ell=1}^{\infty}\big\lfloor m/p^{\ell}\big\rfloor$. Ainsi, on obtient bien 
$$
v_p\left(\mathcal{Q}(n)\right)=v_p\left(\frac{(e_1n)!\cdots(e_{q_1}n)!}{(f_1n)!\cdots(f_{q_2}n)!}\right)=\sum_{\ell=1}^{\infty}\Delta\left(\frac{n}{p^{\ell}}\right)=\sum_{\ell=1}^{\infty}\Delta\left(\left\{\frac{n}{p^{\ell}}\right\}\right).
$$

Nous allons maintenant énoncer un lemme permettant de démontrer \eqref{tre21} que l'on utilisera aussi dans la démonstration du lemme \ref{valuation diff H}.

\begin{lemme}\label{ablanc}
Soit $p$ un nombre premier et $\beta:=\lfloor\log_p(M_{\textbf{e},\textbf{f}})\rfloor$. Pour tout $m\in\mathbb{N}$, $m\geq 1$ et tout $\ell\in\{v_p(m)+1,\dots,v_p(m)+\beta\}$, on a $\{m/p^{\ell}\}\geq 1/M_{\textbf{e},\textbf{f}}$.
\end{lemme}

\begin{proof}
On note $m=p^{v_p(m)}(a+pb)$ avec $a\in\{1,\dots,p-1\}$ et $b\in\mathbb{N}$. On note $b=\sum_{k=0}^{\infty}b_kp^k$, où $b_k\in\{0,\dots,p-1\}$. Pour tout $\ell\in\{v_p(m)+1,\dots,v_p(m)+\beta\}$, on a bien
$$
\left\{\frac{m}{p^{\ell}}\right\}=\left\{\frac{a+bp}{p^{\ell-v_p(m)}}\right\}=\frac{a+p\sum_{k=0}^{\ell-v_p(m)-2}b_kp^k}{p^{\ell-v_p(m)}}\geq\frac{a}{p^{\ell-v_p(m)}}\geq\frac{a}{p^{\beta}}\geq\frac{1}{M_{\textbf{e},\textbf{f}}}.
$$
\end{proof}

\begin{proof}[Démonstration de \eqref{tre21}]
Soit $L\in\{1,\dots,M_{\textbf{e},\textbf{f}}\}$, $\beta:=\lfloor\log_p(M_{\textbf{e},\textbf{f}})\rfloor$, $a\in\{0,\dots,p-1\}$, $K\in\mathbb{N}$ et $j\in\{0,\dots,K\}$ tels que $a\neq 0$ ou $K\neq 0$. On a $v_p(D_1)=p^{\beta}$. On doit donc montrer que
\begin{equation}\label{amontre21}
\mathcal{Q}(K-j)\mathcal{Q}(a+jp)\in p^{\beta}\mathbb{Z}_p.
\end{equation}
Si $K-j=0$, alors $a+jp=a+Kp\neq 0$. Ainsi, on a $K-j\neq 0$ ou $a+jp\neq 0$. Supposons, par exemple, que $K-j\neq 0$. En appliquant le lemme \ref{ablanc} avec $K-j$ à la place de $m$, on obtient que, pour tout $\ell\in\{v_p(K-j)+1,\dots,v_p(K-j)+\beta\}$, on a $\{(K-j)/p^{\ell}\}\geq 1/M_{\textbf{e},\textbf{f}}$. Comme $\Delta$ est positive sur $\mathbb{R}$ et comme, pour tout $x\in[1/M_{\textbf{e},\textbf{f}},1[$, on a $\Delta(x)\geq 1$, on obtient
$$
v_p(\mathcal{Q}(K-j))=\sum_{\ell=1}^{\infty}\Delta\left(\left\{\frac{K-j}{p^{\ell}}\right\}\right)\geq\sum_{\ell=v_p(K-j)+1}^{v_p(K-j)+\beta}\Delta\left(\left\{\frac{K-j}{p^{\ell}}\right\}\right)\geq\beta.
$$
De plus, d'après le critère de Landau, on a $\mathcal{Q}(a+jp)\in\mathbb{N}$ donc on a bien \eqref{amontre21}.

Si $K-j=0$, alors $a+jp\neq 0$ et, de la même manière, le lemme \ref{ablanc} appliqué avec $a+jp$ à la place de $m$ et le fait que $\mathcal{Q}(K-j)\in\mathbb{N}$ donnent bien \eqref{amontre21}.
\end{proof}
Nous allons énoncer un lemme permettant de démontrer le résultat \eqref{lemme Q rho 1} et le lemme \ref{valuation diff H}. Ce lemme est plus fort que le lemme $9$ de \cite{Delaygue} et est la principale amélioration de la preuve du critère pour les applications de type miroir. Ceci nous permet d'obtenir le théorème \ref{RootsL}.

\begin{lemme}\label{lemme 24 généralisé}
Soit $s\in\mathbb{N}$, $s\geq 1$, $a\in\{0,\dots,p^s-1\}$, $M\geq 1$ et $m\in\mathbb{N}$. Soit $L\in\{1,\dots,M\}$ et $\alpha:=\lfloor\log_p(M/L)\rfloor$. Si $\lfloor La/p^s\rfloor\geq 1$ alors, pour tout $u\in\{1,\dots,\lfloor La/p^s\rfloor\}$ et tout $\ell\in\{s,\dots,s+v_p(Lm+u)+\alpha\}$, on a $\left\{(a+mp^s)/p^{\ell}\right\}\geq 1/M$.
\end{lemme}

\begin{proof}
Soit $\ell\in\{s,\dots,s+v_p(Lm+u)\}$. On note $m=\sum_{j=0}^{\infty}m_jp^j$ le développement $p$-adique de $m$. On a 
$$
\left\{\frac{a+mp^s}{p^{\ell}}\right\}=\frac{a+p^s\sum_{j=0}^{\ell-s-1}m_jp^j}{p^{\ell}}.
$$
On a que $p^{\ell-s}$ divise $(u+Lm)$ et donc $p^{\ell-s}$ divise 
$$
u+Lm-L\left(\sum_{j=\ell-s}^{\infty}m_jp^j\right)=u+L\left(\sum_{j=0}^{\ell-s-1}m_jp^j\right).
$$ 
Ainsi, on a
\begin{equation}\label{eqint35}
p^{\ell-s}\leq u+L\left(\sum_{j=0}^{\ell-s-1}m_jp^j\right)\leq\frac{1}{p^s}La+L\left(\sum_{j=0}^{\ell-s-1}m_jp^j\right)=\frac{L}{p^{s-\ell}}\frac{a+p^s\sum_{j=0}^{\ell-s-1}m_jp^j}{p^{\ell}}.
\end{equation}
En multipliant \eqref{eqint35} par $p^{s-\ell}\frac{M}{L}$, on obtient
$$
\frac{M}{L}\leq M\frac{a+p^s\sum_{j=0}^{\ell-s-1}m_jp^j}{p^{\ell}}.
$$
Comme $\alpha\leq\log_p(M/L)$, pour tout $i\in\{0,\dots,\alpha\}$, on a $p^i\leq M/L$ et donc
$$
p^i\leq M\frac{a+p^s\sum_{j=0}^{\ell-s-1}m_jp^j}{p^{\ell}}.
$$
Ainsi, pour tout $\ell\in\{s,\dots,s+v_p(Lm+u)\}$ et tout $i\in\{0,\dots,\alpha\}$, on a
$$
1\leq M\frac{a+p^s\sum_{j=0}^{\ell-s-1}m_jp^j}{p^{\ell+i}}\leq M\frac{a+p^s\sum_{j=0}^{\ell+i-s-1}m_jp^j}{p^{\ell+i}}=M\left\{\frac{a+mp^s}{p^{\ell+i}}\right\}
$$
et on obtient bien que, pour tout $\ell\in\{s,\dots,s+v_p(Lm+u)+\alpha\}$, on a $\left\{(a+mp^s)/p^{\ell}\right\}\geq 1/M$.
\end{proof}

Nous allons maintenant appliquer le lemme \ref{lemme 24 généralisé} pour démontrer \eqref{lemme Q rho 1}.

\begin{proof}[Démonstration de \eqref{lemme Q rho 1}]
Soit $L\in\{1,\dots,M_{\textbf{e},\textbf{f}}\}$, $\alpha:=\lfloor\log_p(M_{\textbf{e},\textbf{f}}/L)\rfloor$, $a\in\{0,\dots,p-1\}$ et $j\in\mathbb{N}$. Comme $v_p(D_L)=\alpha$, il faut montrer que 
$$
\mathcal{Q}(a+jp)\sum_{i=1}^{\lfloor La/p\rfloor}\frac{1}{Lj+i}\in p^{\alpha+1}\mathbb{Z}_p.
$$ 
Si $\lfloor La/p\rfloor=0$, c'est évident. Supposons que $\lfloor La/p\rfloor\geq 1$. En appliquant le lemme \ref{lemme 24 généralisé} avec $s=1$, $m=j$ et $M=M_{\textbf{e},\textbf{f}}$, on obtient que, pour tout $i\in\{1,\dots,\lfloor La/p\rfloor\}$ et tout $\ell\in\{1,\dots,1+v_p(i+Lj)+\alpha\}$, on a $\{(a+jp)/p^{\ell}\}\geq 1/M_{\textbf{e},\textbf{f}}$ et donc $\Delta(\{(a+jp)/p^{\ell}\})\geq 1$. Comme $\Delta$ est positive sur $\mathbb{R}$, cela donne
$$
v_p(\mathcal{Q}(a+jp))=\sum_{\ell=1}^{\infty}\Delta\left(\left\{\frac{a+jp}{p^{\ell}}\right\}\right)\geq\sum_{\ell=1}^{1+v_p(Lj+i)+\alpha}\Delta\left(\left\{\frac{a+jp}{p^{\ell}}\right\}\right)\geq 1+v_p(Lj+i)+\alpha,
$$
ce qui achève la preuve de \eqref{lemme Q rho 1}.
\end{proof}

\begin{proof}[Démonstration du lemme \ref{valuation diff H}]
Soit $L\in\{1,\dots,M_{\textbf{e},\textbf{f}}\}$, $\alpha:=\lfloor\log_p(M_{\textbf{e},\textbf{f}}/L)\rfloor$ et $m$ et $s$ dans $\mathbb{N}$. Il faut montrer que 
$$
p^{s+1}g_p(m)(H_{Lmp^s}-H_{L\lfloor m/p\rfloor p^{s+1}})\in p^{\alpha+1}\mathbb{Z}_p.
$$ 

Si $m=0$, alors c'est évident. On suppose, dans la suite de la démonstration, que $m\geq 1$. On écrit $m=b+qp$, où $b\in\{0,\dots,p-1\}$ et $q\in\mathbb{N}$. On a alors $Lmp^s=Lbp^s+Lqp^{s+1}$ et $L\lfloor m/p\rfloor p^{s+1}=Lqp^{s+1}$. Ainsi, on obtient
$$
H_{Lmp^s}-H_{L\lfloor m/p\rfloor p^{s+1}}=\sum_{j=1}^{Lbp^s}\frac{1}{Lqp^{s+1}+j}\equiv\sum_{i=1}^{\lfloor Lb/p\rfloor}\frac{1}{Lqp^{s+1}+ip^{s+1}}\mod \frac{1}{p^s}\mathbb{Z}_p
$$
et donc 
$$
p^{s+1}g_p(m)(H_{Lmp^s}-H_{L\lfloor m/p\rfloor p^{s+1}})\equiv g_p(b+qp)\sum_{i=1}^{\lfloor Lb/p\rfloor}\frac{1}{Lq+i}\mod pg_p(m)\mathbb{Z}_p.
$$
On note $\beta:=\lfloor\log_p(M_{\textbf{e},\textbf{f}})\rfloor$. D'après le lemme \ref{ablanc}, pour tout $\ell\in\{v_p(m)+1,\dots,v_p(m)+\beta\}$, on a $\{m/p^{\ell}\}\geq 1/M_{\textbf{e},\textbf{f}}$. On obtient
$$
v_p(g_p(m))=\sum_{\ell=1}^{\infty}\textbf{1}_{[1/M_{\textbf{e},\textbf{f}},1[}\left(\left\{\frac{m}{p^{\ell}}\right\}\right)\geq\sum_{\ell=v_p(m)+1}^{v_p(m)+\beta}\textbf{1}_{[1/M_{\textbf{e},\textbf{f}},1[}\left(\left\{\frac{m}{p^{\ell}}\right\}\right)\geq\beta\geq\alpha.
$$ 
Ainsi, on a
$$
p^{s+1}g_p(m)(H_{Lmp^s}-H_{L\lfloor m/p\rfloor p^{s+1}})\equiv g_p(b+qp)\sum_{i=1}^{\lfloor Lb/p\rfloor}\frac{1}{Lq+i}\mod p^{\alpha+1}\mathbb{Z}_p
$$
et il ne nous reste plus qu'à montrer que 
$$
g_p(b+qp)\sum_{i=1}^{\lfloor Lb/p\rfloor}\frac{1}{Lq+i}\in p^{\alpha+1}\mathbb{Z}_p.
$$ 
Si $\lfloor Lb/p\rfloor=0$, c'est évident. Supposons que $\lfloor Lb/p\rfloor\geq 1$. En appliquant le lemme \ref{lemme 24 généralisé} avec $s=1$, $M=M_{\textbf{e},\textbf{f}}$, $a=b$ et $q$ à la place de $m$, on obtient que, pour tout $i\in\{1,\dots,\lfloor Lb/p\rfloor\}$ et tout $\ell\in\{1,\dots,1+v_p(i+Lq)+\alpha\}$, on a $\{(b+qp)/p^{\ell}\}\geq1/M_{\textbf{e},\textbf{f}}$ et donc
\begin{multline*}
v_p(g_p(b+qp))=\sum_{\ell=1}^{\infty}\textbf{1}_{[1/M_{\textbf{e},\textbf{f}},1[}\left(\left\{\frac{b+qp}{p^{\ell}}\right\}\right)\\
\geq\sum_{\ell=1}^{1+v_p(Lq+i)+\alpha}\textbf{1}_{[1/M_{\textbf{e},\textbf{f}},1[}\left(\left\{\frac{b+qp}{p^{\ell}}\right\}\right)\geq 1+v_p(Lq+i)+\alpha,
\end{multline*}
ce qui achève la preuve du lemme.
\end{proof}

\section{Démonstration du théorème \ref{cor1} et du corollaire \ref{cor2}}\label{partie4}

\subsection{Démonstration du théorème \ref{cor1}}

Soit $\textbf{e}$ et $\textbf{f}$ deux suites d'entiers strictement positifs disjointes vérifiant $|\textbf{e}|=|\textbf{f}\,|$ et telles que, pour tout $x\in[1/M_{\textbf{e},\textbf{f}},1[$, on a $\Delta_{\textbf{e},\textbf{f}}(x)\geq 1$. Soit $\theta\geq 1$ un diviseur de $M_{\textbf{e},\textbf{f}}$ tel que, pour tout élément $L$ de $\textbf{e}$ et $\textbf{f}$, on a que $\theta/\textup{pgcd}(L,\theta)$ divise $D_L$. On doit montrer que
$$
(z^{-1}q_{\textbf{e},\textbf{f}}(z))^{1/\theta}\in\mathbb{Z}[[z]].
$$
En notant $\textbf{e}=(e_1,\dots,e_{q_1})$ et $\textbf{f}=(f_1,\dots,f_{q_2})$, on rappelle qu'on a
$$
z^{-1}q_{\textbf{e},\textbf{f}}(z)=\frac{\prod_{i=1}^{q_1}q_{e_i,\textbf{e},\textbf{f}}(z)^{e_i}}{\prod_{i=1}^{q_2}q_{f_i,\textbf{e},\textbf{f}}(z)^{f_i}}.
$$
Il suffit donc de montrer que, pour tout $L\in\{e_1,\dots,e_{q_1},f_1,\dots,f_{q_2}\}$, on a 
$$
q_{L,\textbf{e},\textbf{f}}(z)^{L/\theta}=(q_{L,\textbf{e},\textbf{f}}(z)^{\textup{pgcd}(L,\theta)/\theta})^{L/\textup{pgcd}(L,\theta)}\in\mathbb{Z}[[z]],
$$
ce qui est impliqué, comme $L/\textup{pgcd}(L,\theta)\in\mathbb{N}$, par le fait que, pour tout élément $L$ de $\textbf{e}$ et $\textbf{f}$, on a $q_{L,\textbf{e},\textbf{f}}(z)^{\textup{pgcd}(L,\theta)/\theta}\in\mathbb{Z}[[z]]$. 

Soit $L\in\{e_1,\dots,e_{q_1},f_1,\dots,f_{q_2}\}$ et $k_L\in\mathbb{N}$ tel que $D_L=\frac{\theta}{\textup{pgcd}(L,\theta)}k_L$. D'après le théorème~\ref{RootsL}, on a $q_{L,\textbf{e},\textbf{f}}(z)^{1/D_L}\in\mathbb{Z}[[z]]$ donc on a bien
$$
q_{L,\textbf{e},\textbf{f}}(z)^{\textup{pgcd}(L,\theta)/\theta}=(q_{L,\textbf{e},\textbf{f}}(z)^{1/D_L})^{k_L}\in\mathbb{Z}[[z]].
$$

\subsection{Démonstration du corollaire 1}

Soit $\textbf{e}$ et $\textbf{f}$ deux suites d'entiers strictement positifs disjointes vérifiant $|\textbf{e}|=|\textbf{f}\,|$ telles que, pour tout $x\in[1/M_{\textbf{e},\textbf{f}},1[$, on a $\Delta_{\textbf{e},\textbf{f}}(x)\geq 1$ et telles que tout élément de $\textbf{e}$ et $\textbf{f}$ divise $M_{\textbf{e},\textbf{f}}$. Nous allons appliquer le théorème \ref{cor1} avec $\theta=M_{\textbf{e},\textbf{f}}$. Pour tout élément $L$ de $\textbf{e}$ et $\textbf{f}$, on a $\lfloor M_{\textbf{e},\textbf{f}}/L\rfloor=M_{\textbf{e},\textbf{f}}/L$ donc $M_{\textbf{e},\textbf{f}}/L$ divise $D_L$. Ainsi $\theta/\textup{pgcd}(L,\theta)=M_{\textbf{e},\textbf{f}}/\textup{pgcd}(L,M_{\textbf{e},\textbf{f}})=M_{\textbf{e},\textbf{f}}/L$ divise $D_L$ et on a bien
\begin{equation}\label{resulti}
(z^{-1}q_{\textbf{e},\textbf{f}}(z))^{1/M_{\textbf{e},\textbf{f}}}\in\mathbb{Z}[[z]].
\end{equation}

Démontrons maintenant la conjecture de Zhou. Soient $k_1,\dots,k_n$ des entiers strictement positifs vérifiant $\frac{1}{k_1}+\cdots+\frac{1}{k_n}=1$. Soit $k$ le plus petit multiple commun de $k_1,\dots,k_n$ et, pour tout $i\in\{1,\dots,n\}$, $w_i:=k/k_i$. En notant $\textbf{e}:=(k)$ et $\textbf{f}:=(w_1,\dots,w_n)$, on obtient que $|\textbf{e}|-|\textbf{f}\,|=k-\sum_{i=1}^nw_i=0$. De plus, Zhou montre dans la partie $2$ de \cite{Zhou1} que, pour tout $x\in[1/k,1[$, on a $\Delta_{\textbf{e},\textbf{f}}(x)\geq 1$. Enfin, pour tout $i\in\{1,\dots,n\}$, $w_i$ divise $k$. On peut donc appliquer \eqref{resulti} qui donne bien \eqref{corZhou} puisque $M_{\textbf{e},\textbf{f}}=k$. Ce qui termine la preuve du corollaire \ref{cor2}.

\address{E. Delaygue, Institut Fourier, CNRS et Université Grenoble 1, 100 rue des Maths, BP 74, 38402 Saint-Martin-d'Hères cedex, France. Email : Eric.Delaygue@ujf-grenoble.fr}

\end{document}